\documentclass[11pt]{article}

\usepackage[utf8]{inputenc}  
\usepackage[T1]{fontenc}

\usepackage{amsfonts,amsmath,amssymb,latexsym}
\usepackage{bbm} 
\usepackage[dvips]{graphicx,epsfig}

\def\real{\mathbb{R}}

\def\mat1{\mathbb{I}}
\def\bmat1{\boldsymbol{\mathbb{I}}}

\def\dzero{\mbox{{\scriptsize $\mathbb{O}$}}}  
\def\dun{\mathbbm 1}

\def\limk{\lim_{k \rightarrow +\infty}}

\newcommand { \Iletter}[1] {I\kern-0.10em #1 }

\def\bit{\begin{itemize}}
\def\eit{\end{itemize}}
\def\ben{\begin{enumerate}}
\def\een{\end{enumerate}}
\def\bde{\begin{description}}
\def\ede{\end{description}}

\def\bar{\begin{array}}
\def\ear{\end{array}}
\def\beq{\begin{equation}}
\def\eeq{\end{equation}}
\def\bfi{\begin{figure}[hbt] \begin{center}}
\def\efi{\end{center} \end{figure}}
\def\noi{\noindent}
\def\bce{\begin{center}}
\def\ece{\end{center}}
\newcommand{\proof}{{\bf Preuve. }}
\newcommand{\cqfd}{\hfill $\Box$}

\newtheorem {proc} {Mode de remplissage pour obtenir une forme matricielle}[section]

\newtheorem {theo} {Théorème}[section]

\newtheorem {lemm} {Lemme}[section]

\newtheorem {propo} {Proposition}

\newtheorem {defi} {Définition}

\newtheorem {rema} {Remarque} [section]

\begin{document}


\title{Looking for all solutions of the Max Atom Problem (MAP)}
\author{L. Truffet \\
  IMTA/DAPI \\
@: laurent.truffet@imt-atlantique.fr; ltruff@proton.me}

\maketitle

\begin{abstract}
This present paper provides the absolutely necessary corrections to the previous
  work entitled {\it  A polynomial Time Algorithm to Solve The
Max-atom Problem} (arXiv:2106.08854v1).

The max-atom-problem (MAP)  deals with system of
  scalar inequalities (called atoms or max-atom) of the form: $x \leq a + \max(y,z)$. Where
  $a$ is a real number and $x,y$ and $z$ belong to the set 
  of the variables of the whole MAP. A max-atom is said to be positive if
  its scalar $a$ is $\geq 0$ and stricly negative if its
  scalar $a <0$. A MAP will be said to be positive
  if all atoms are positive. In the case of non positive MAP
  we present a saturation principle for
  system of vectorial inequalities of the form
  $x \leq A x + b$ in the so-called $(\max,+)$-algebra assuming
  some properties on the matrix $A$. Then, we apply such
  principle to explore all non-trivial solutions (ie $\neq -\infty$).
  We deduce a strongly polynomial method to express all
  solutions of a non positive MAP. In the case a positive MAP which has
  always the vector $x^{1}=(0)$ as trivial solution we show that looking
  for all solutions requires the enumeration of all elementary circuits
  in a graph associated with the MAP. However, we propose a strongly
  polynomial method wich provides some non trivial solutions.

\end{abstract}

\noi
 Non positive MAP: strongly polynomial method.
    Positive MAP: elementary circuits; monomial matrix; non-trivial solutions.

\section*{Remerciements}

Un grand merci à Chams Lahlou pour son énergie passée
et qui a grandement contribué à l'exploration puis l'élagage de
bon nombre de fausses bonnes idées. \\

Un grand merci à James Ledoux qui, malgré un emploi du temps de
dingue, a relu les passages importants de ce papier. \\

Un grand merci aussi à Tom Van Dijk pour son travail qui a permis de pointer les erreurs contenues
dans l'ancien travail intitulé {\it  A polynomial Time Algorithm to Solve The
  Max-atom Problem} (arXiv:2106.08854v1). \\

Un immense merci à Nadjib pour la lumière qu'il m'a apporté,  moi
qui pensais que chaque jour était un miracle. Mais non!
{\em Chaque jour est un don}. 

\tableofcontents

\section{Introduction}
\label{secIntro}
Le problème du max-atom ou MAP a une formulation très simple qui est 
la suivante. Il faut trouver l'ensembe des solutions d'un système 
de max-atomes (ou atomes tout simplement), c'est-à-dire des inégalités du type: 
$x \leq a + \max(y,z)$, ou bien seulement $x \leq a + y$, avec $a$ qui 
est un réel et 
les variables sont $x,y,z$. Dans le cas général il peut y avoir un 
nombre $n$ de variables et un nombre $m$ d'atomes.

\noi{\bf Un pré-traitement élémentaire}. De ce système
d'inégalités, on aura pris soin au préalable
de retirer les max-atoms du type:
\[
x \leq a + \max(x,y)
\]
avec $a \geq 0$, car ce type d'inégalités a la valeur logique 'vraie'.

Un max-atom du type $x \leq a + \max(x,y)=\max(a+x, a+y)$ sera transformer
en l'inégalité: $x \leq a + y$ lorsque $a < 0$. En effet, dans ce cas
l'inégalité $x \leq a+x$ n'est valable que si $x=-\infty$. Or,
la solution $(-\infty)$ est toujours solution d'un système d'atomes. 

\subsection{Exemple de base}
\label{subEDB}

Considérons quatre variables: $x_{1}, \ldots, x_{4}$ et le 
système $S$ d'atomes suivants:

\[
S=\{x_{3} \leq (-10) + x_{1}, x_{4} \leq (-1) + \max(x_{2},x_{3}), x_{2} \leq x_{4}, 
x_{4} \leq 25 + \max(x_{2},x_{3}) \}.
\]

Ce MAP de base, ie le système $S$ est empreinté à
\cite{kn:Bezemetal2008}, Example $1$ avec
les correspondances entre variables suivantes:

\[
x_{1} \leftrightarrow u, \; x_{2} \leftrightarrow y, \; x_{3} \leftrightarrow x, \; x_{4} \leftrightarrow z.
\]

Dans la suite on pourra modifier $S$ comme suit. Soit en considérant
le système $S'=S \cup \{x_{1} \leq 9 + \max(x_{2},x_{3})\}$ (cf.\cite{kn:Bezemetal2008}, Example 2), soit en considérant
le système $S''= S \cup \{x_{1} \leq 12 +\max(x_{2},x_{3})\}$.

\section{Algèbre $(\max,+)$: rappels succints et notations}
Pour de plus amples précisions nous renvoyons le lecteur
à e.g \cite{kn:Bac-cooq}.

L'algèbre $(\max,+)$ désigne l'ensemble de référence, noté
$\real_{\dzero}:= \real \cup \{-\infty \}$ associé aux
opérations de seront: $(a,b) \mapsto \max(a,b)$ qui jouera le rôle de
``l'addition'', et $(a,b) \mapsto a + b$ qui jouera le rôle de la multiplication.

Dans la suite comme nous resterons dans cette algèbre nous emploierons les notations
suivantes:
\bit
\item $a + b := \max(a,b)$ (en opération classique),
\item $a \cdot b   := a+ b$ (en opération classique). $a \cdot b$ sera
  raccourcie en $ab$ tout simplement.
  \eit
  
  $\dzero :=-\infty$ est l'élément neutre de $+$, $\dun :=0$ est l'élément neutre
  de $\cdot$.

  La notation puissance sera utilisée et aura la signification suivante:

  \[
  a^{b}:= a \times b \mbox{ (en opération classique)}.
  \]

  $\real_{\dzero}$ sera muni de l'ordre naturel $\leq$ défini par:
  \[
a \leq b \Leftrightarrow a + b = b. 
  \]

  Nous aurons également besoin de l'opérateur $\min$ noté dans la suite $\wedge$.

  Les opérations et relations précédentes se généralisent pour les matrices comme
  suit:
  \bit
\item addition: $[a_{i,j}] + [b_{i,j}] := [a_{i,j} + b_{i,j}]$ (rappel: $a_{i,j}+b_{i,j}
  = \max(a_{i,j},b_{i,j}$).
\item minimum: $[a_{i,j}] \wedge [b_{i,j}] := [a_{i,j} \wedge b_{i,j}]$
\item produit: $C = A B$ a pour entrée $(i,j)$:
  $c_{i,j}:= \Sigma_{k} a_{i,k} b_{k,j}$. Cette
  expression correspond à $\max_{k}(a_{i,k}+b_{k,j})$.
\item comparaison $A \leq B$ équivaut à $\forall i,j: a_{i,j} \leq b_{i,j}$

  \eit
  Dans cette algèbre, $I_{n}$ désignera la $n \times n$-matrice identité
  (tous ses termes hors diagonale sont $\dzero$, et ses termes diagonaux
  sont $\dun$). Et la
  matrice $O_{m,n}$ désignera la $m \times n$-matrice nulle (tous ses termes
  sont $\dzero$).

  Nous utiliserons également la résiduation. Cette théorie est le pendant de l'inversion matricielle
  en algèbre linéaire. Soient $m,n,q$ sont des entiers $\geq 1$. Considérons
  les matrices données $A$: $m \times n$ et $B$: $m \times q$. Il s'agit
  de trouver la plus grande  $n \times q$-matrice $X$
  (au sens de l'ordre partiel sur les matrices, $\leq$) qui vérifie
  l'inégalité matricielle suivante:

  \begin{equation}
    \label{eqpourlaresiduation}
AX \leq B.
  \end{equation}
  Cette inégalité se résoud facilement en utilisant le fait que:
  $a + b \leq c \Leftrightarrow (a \leq c) \mbox{ et } (b\leq c)$. Et
  l'on peut écrire:

  \begin{equation}
    \label{eqmatresiduee}
X \leq A \backslash B.
  \end{equation}
  L'opérateur $A \backslash (\cdot)$ joue le rôle de $A^{-1}(\cdot)$ en
  l'algèbre linéaire. La matrice $A \backslash =[(a\backslash)_{i,j}]$ est une $n \times m$-matrice
  telle que

  \begin{equation}
    \label{eqAbackslash}
\forall i=1, \ldots,n, \forall j=1, \ldots, m, (a\backslash)_{i,j} := a^{-1}_{j,i},
    \end{equation}
  avec la convention: $\dzero^{-1}:= +\infty$. L'application bijective $ \cdot^{-1}$ a pour ensemble de
  départ $\real \cup \{-\infty \}$ et pour ensemble d'arrivée $\real \cup \{+\infty \}$.
  La $n \times q$-matrice $A \backslash B$ se définit comme un produit de $A\backslash$ avec $B$, mais
  dans l'algèbre $(\real \cup \{+\infty \}, \min, +)$ où $\min$ est considéré comme une addition et
  $+$ joue toujours le rôle de la multiplication. Cela signifie que l'entrée $(i,j)$ de
  $A \backslash B$ notée $(A \backslash B)_{i,j}$ est définie par:

  \begin{equation}
    \label{eqAbackslashBij}
(A \backslash B)_{i,j} := \min_{k=1, \ldots, m}((a\backslash)_{i,k} + b_{k,i}), 
    \end{equation}
pour $i=1, \ldots, n$ et $j=1, \ldots, q$.
  
  \section{Lemme de saturation ou maximisation}
  \label{secLemmSat}

  Dans cette section nous étudions les solution en $x$ de l'inégalité vecorielle
  suivante suivante:
  
  \begin{equation}
    \label{masterineq}
x \leq A x + b
  \end{equation}
  
  où $A$ est une $n \times n$-matrice, $b$ vecteur colonne à $n$ composantes
  et $x$ vecteur colonne à $n$ inconnues.
  
  \begin{lemm}[Saturation]
    \label{lemSAT}
    Si $\exists \limk A^{k} = O$ (matrice nulle) alors la plus grande
    solution de l'inéquation (\ref{masterineq}) est:
    \[
x= A^{*} b
\]
avec $A^{*}:= I + A + A^{2}+ A^{3} + \cdots = I + A + A^{2}+ \cdots + A^{n-1}$.
Cette plus grande solution coïncide avec la solution de
l'équation $x= Ax + b$, qui correspond au cas de
la saturation ou maximisation du système d'inégalités (\ref{masterineq}).
    \end{lemm}
  \proof
  Cette preuve est fondée sur le fait que la fonction $x \mapsto A x$ est
  croissante. Elle se fait par récurrence comme dans le cas du système
  d'égalités. Nous la reproduisons ici.

  Si $x \leq A x + b$ (1) alors comme $x \mapsto Ax$ est croissante
  on $A x \leq A (A x + b)= A^{2} x + A b$ (1'). Comme $+=\max$ est
  également croissante on peut substituer $A x$ dans (1) par sa borne
  supérieure: $A^{2} x + A b$. On obtient alors:
  $ x \leq A^{2} x + (I + A) b$ (2). Comme la fonction $x \mapsto A x$
  est croissante, à partir de (1') on a (multiplication par $A$ à gauche
  et à droite du signe $\leq$): $A^{2} x \leq A^{3} x + A^{2} b$ (2').
  Comme $+ = \max$ est croissante, on peut remplacer $A^{2} x$ par sa borne
  supérieure issue de (2'). Il vient $x \leq A^{3} x +(I + A + A^{2}) b$. Et ainsi de suite, au rang $k$ on pourra écrire:
  \[
x \leq A^{k} x + (I + A + \cdots + A^{k-1}) b. 
  \]

  Par passage à la limite lorsque $k \rightarrow \infty$ on a:

  \[
x \leq A^{*} b.
  \]
  Et il est trivial de vérifier que cette borne supérieure $ A^{*} b$
  est solution de $x= A x + b$.
  
  \cqfd

  Examinons maintenant sous quelles conditions nous avons:

  \[
\exists \limk A^{k}= O.
  \]
  A toute $n \times n$-matrice $A$ on peut lui associer un graphe valué 
  $G(A)$ comme suit. Les sommets du graphe $G(A)$ seront les
  lignes $\{1, \ldots, n\}$. Il existera un arc $(j,i)$ si et seulement
  si $a_{i,j} \neq \dzero$. Et sa valuation sera dans ce cas $a_{i,j}$.
  Un chemin de longueur $k \geq 1$ (ie une suite de $k$ arcs consécutifs)
  $i_{0} \rightarrow i_{1} \rightarrow
  \cdots \rightarrow i_{k}$ aura pour
  poids: $a_{i_{1},i_{0}} a_{i_{2},i_{1}} \ldots a_{i_{k}, i_{k-1}}$. Dans le cas
  particulier où $i_{k}=i_{0}$ ce chemin sera appelé circuit. Si de plus les
  sommets $i_{0}, i_{1}, \ldots, i_{k}$ sont tous distincts alors
  le circuit est dit élémentaire.
  Le terme noté $a_{i,j}^{(k)}$ correspondant à
  l'entrée $(i,j)$ de la matrice $A^{k}$ est le maximum des poids de tous
  les chemins de longueur $k$ (exactement) allant de $j$ vers $i$.

  Nous avons les cas.

  \bit
\item Cas $1$. La matrice $A$ est nilpotente. De manière équivalente nous pouvons dire que
  $G(A)$ est un graphe orienté sans circuits ou DAG.
  
\item Cas $2$. La matrice $A$ est irréductible. Ou de manière équivalente
  son graphe $G(A)$ est fortement connexe. Nous avons le
  résultat classique suivant:

  \[
  \mbox{ tous les circuits de $G(A)$ sont de poids $< \dun$} \Rightarrow
  \limk A^{k}= O.
  \]
\item Cas $3$. La matrice $A$ est dite réductible. Dans ce cas elle
  admet une décomposition par blocs quitte à permuter les lignes et les
  colonnes de la forme bloc-triangulaire suivante:
  \[
  A= \left(\bar{ccccc} A_{1,1} & O & \cdots & \cdots & O \\
  A_{2,1} & A_{2,2} & O & \cdots & O \\
  \vdots & \vdots & \ddots & \ddots & \vdots \\
  \vdots & \vdots & \vdots & \ddots & \vdots \\
  A_{r,1} & \cdots & \cdots & A_{r,r-1} & A_{r,r}
  \ear \right).
  \]

  Avec $r \geq 2$ et $A_{1,1}, \ldots, A_{r,r}$ matrices carrées irréductibles
  ou nulles. 

Nous avons alors le résultat classique suivant:
\[
\forall i=1, \ldots , r \; \exists \limk A_{i,i}^{k}=O \Rightarrow \limk A^{k}=O.
\]

\item Cas $4$. Toutes les entrées de la matrice $A$ sont
  $ < \dun$.
\eit

\section{Mise sous forme matricielle du MAP}
\label{secMAPmatrice}
Nous avons vu que le problème MAP consiste en un ensemble de
conjonctions d'atomes à résoudre. Nous considèrerons que nous
avons un problème avec $n$ variables $x=(x_{i})_{i=1}^{n}$ et
$m$ atomes. Nous allons établir une procédure qui permet de
construire $L \geq 1$ $n \times n$-matrices $A_{1}, \ldots, A_{L}$
telles que le MAP est équivalent au système d'inégalités:
\[
x \leq A_{1} x \wedge \ldots \wedge A_{L} x. 
\]
Comme le 'et logique' est associatif, idempotent et commutatif
nous pourront imposer certaines conditions sur la façon de remplir ces
matrices. 

\subsection{Formulation algébrique de l'exemple de base}
Considérons le système d'atomes $S$. Nous remarquons que la variable
$x_{1}$ n'apparait jamais à gauche du signe $\leq$. Nous rajoutons
alors l'atome trivial $x_{1} \leq x_{1}$. Le nouveau système considéré
est alors $S \cup \{x_{1} \leq x_{1}\}$. Ce système s'écrit alors:

\[
\bar{ccc}
x_{1} & \leq & x_{1} \\
x_{2} & \leq & x_{4} \\
x_{3} & \leq & (-10) + x_{1} \\
x_{4} & \leq & \min((-1) + \max(x_{2}, x_{3}), 25+ \max(x_{2}, x_{3}))
\ear
\]
Cette forme est très peut lisible. Nous passons aux notations
maxplusiennes, il vient:
\[
\bar{ccc}
x_{1} & \leq & x_{1} \\
x_{2} & \leq & x_{4} \\
x_{3} & \leq & 10^{-1} x_{1} \\
x_{4} & \leq & (1^{-1} (x_{2}+x_{3})) \wedge (25 (x_{2}+x_{3})).
\ear
\]
Ici nous pourrions utiliser la simplification
$(1^{-1} (x_{2}+x_{3})) \wedge (25 (x_{2}+x_{3})) = 1^{-1} (x_{2}+x_{3})$
mais nous ne le ferons pas car c'est un cas particulier. Comme
$x_{4}$ doit être inférieure à deux fonctions, nous pouvons
utiliser le fait que le 'et logique' est idempotent pour avoir
la formulation suivante:
\[
\bar{ccc}
x_{1} & \leq & x_{1} \wedge x_{1} \\
x_{2} & \leq & x_{4} \wedge x_{4} \\
x_{3} & \leq & 10^{-1} x_{1} \wedge 10^{-1} x_{1}\\
x_{4} & \leq & (1^{-1} (x_{2}+x_{3})) \wedge (25 (x_{2}+x_{3})).
\ear
\]
Pour chaque fonction à droite de $\leq$ nous pouvons
lui associer une reprsésentation vectorielle en introduisant des
vecteurs lignes. Ainsi, la fonction $x_{1} = a_{1} x$ avec
$a_{1}=(\dun, \dzero, \dzero, \dzero)$. La fonction $x_{4}= a_{2} x$,
avec $a_{2}=(\dzero, \dzero, \dzero , \dun)$. La fonction $10^{-1} x_{1}=
a_{3} x$ avec $a_{3} = (10^{-1}, \dzero, \dzero, \dzero)$. La fonction
$1^{-1} (x_{2}+x_{3}= a_{4}^{1} x$, avec $a_{4}^{1}=(\dzero, 1^{-1}, 1^{-1}, \dzero)$.
Et la fonction $25 (x_{2}+x_{3})= a_{4}^{2} x$, avec $a_{4}^{2}=(\dzero, 25, 25, \dzero)$.

Le système d'inégalités précédent se réécrit alors:

\[
\bar{ccc}
x_{1} & \leq & a_{1} x \wedge a_{1} x \\
x_{2} & \leq & a_{2} x \wedge a_{2} x \\
x_{3} & \leq & a_{3} x \wedge a_{3} x \\
x_{4} & \leq & a_{4}^{1} x \wedge a_{4}^{2} x.
\ear
\]

Si l'on prend par exemple
$A_{1}=\left(\bar{c} a_{1} \\ a_{2} \\ a_{3} \\ a_{4}^{1} \ear \right)$ et
$A_{2} =\left(\bar{c} a_{1} \\ a_{2} \\ a_{3} \\ a_{4}^{2} \ear \right)$
alors le système précédent se réécrit maintenant:

\[
x \leq A_{1} x \wedge A_{2} x.
\]
Le problème est que cette représentation n'est pas unique. En effet,
on pourrait également prendre:
$A_{1}=\left(\bar{c} a_{1} \\ a_{2} \\ a_{3} \\ a_{4}^{2} \ear \right)$ et
$A_{2}=\left(\bar{c} a_{1} \\ e_{2} \\ e_{3} \\ a_{4}^{1} \ear \right)$, avec
$e_{2}=(\dzero, \dun, \dzero, \dzero)$ et $e_{3}=(\dzero, \dzero, \dun, \dzero)$
qui désigne le deuxième et le troisième vecteurs de la base canonique
de $\real_{\dzero}^{4}$.

En utilisant les conditions matricielles pour
pouvoir appliquer le lemme de saturation nous allons montrer dans le
cas général (sous-section suivante) comment choisir
la représentation matricielle que nous appellerons la plus
pénalisante.

\subsection{Cas général}
\label{subCasGenAlgebric}
Nous considérons un ensemble $\mathcal{A}$ de $m$ atomes utilisant
$n$ variables stockées dans le vecteur colonne $x=(x_{i})_{i=1}^{n}$.

  Les vecteurs lignes de la matrice identité $I_{n}$ seront
  notés $e_{1}, \ldots, e_{n}$. Ce sont aussi les vecteurs
  de la base canonique de $\real_{\dzero}^{n}$. A tout atome
  $x_{i} \leq a (x_{j} + x_{k})$ de $\mathcal{A}$ on peut
  lui associer le vecteur ligne $a_{i}:= a e_{j} + a e_{k}$ (l'indice
  $i$ de $a_{i}$ pour rappeler la variable $x_{i}$). En effet, comme
  $a_{i} x= a x_{j} + a x_{k}$, l'atome  $x_{i} \leq a (x_{j} + x_{k})$
  se réécrit: $x_{i} \leq a_{i} x$. 
  Cet atome sera dit strictement négatif si $a < \dun$ (ce qui veut
  dire aussi que toutes les composantes de $a_{i}$ sont $ < \dun$).
  Et il sera dit positif si $a \geq \dun$. A
  chaque variable $x_{i}$ on lui associe:

  \bit
\item sa liste d'atomes strictement négatifs $\mathcal{L}^{-}(x_{i})$ à laquelle
  on rajoute à la fin de la liste l'atome $\{x_{i} \leq x_{i}\}$,
  \item sa liste d'atomes positifs  $\mathcal{L}^{+}(x_{i})$ à laquelle
  on rajoute en dernier l'atome $\{x_{i} \leq x_{i}\}$.
  \eit 
  Dans les faits on stockera les vecteurs lignes associés aux atomes. A titre
  d'exemple l'atome
  $\{x_{i} \leq x_{i}\}$ sera associé au vecteur ligne $e_{i}$.

  Ceci nous amène à la procédure de remplissage des matrices
  suivante. La forme matricielle se justifie par le fait que le 'et logique'
  est associatif. 

  \noi
  \begin{proc}
    \label{procRempli}
   Cette procédure s'effectue comme suit.
   
 \bit
\item Remplissage de la partie dite strictement négative du MAP. \\
  
 A l'étape $1$, l'inégalité $x \leq A_{1} x$ s'obtient en construisant
 la matrice $A_{1}$ de la manière suivante:

  \[
  \forall i=1, \ldots, n, \; a_{i}^{1}= \mbox{ premier vecteur de la liste
   $\mathcal{L}^{-}(x_{i})$}
  \]
  On ne retire de la liste $\mathcal{L}^{-}(x_{i})$ que les vecteurs $\neq e_{i}$.
  On obtient une nouvelle liste encore notée abusivement $\mathcal{L}^{-}(x_{i})$.
  
  Tant que $\exists i = 1, \ldots,n: \mathcal{L}^{-}(x_{i}) \neq \{e_{i}\}$,
  on construit des matrices suivant la même procédure.

  Dans la suite nous admettrons que l'on peut construire $l$ matrices
  $A_{1}, \ldots, A_{l}$, $l \geq 1$.
\item Construction de la partie positive du MAP
  
  De même on construit la liste des matrices $A_{l+1}, \ldots, A_{L}$
  avec $L \geq l+1$ en utilisant la procédure précédente mais en
  remplaçant $\mathcal{L}^{-}(x_{i})$ par $\mathcal{L}^{+}(x_{i})$.
  \eit
  
  Au final, nous pouvons écrire que le MAP est équivalent à
  $x \leq A_{1} x$ et ...et $x \leq A_{l} x$ et $x \leq A_{l+1} x$ et ...
  et $x \leq A_{L} x$, ce qui revient à écrire de manière plus condensée:

  \begin{equation}
x \leq A_{1} x \wedge \ldots \wedge A_{l}x \wedge A_{l+1}x \wedge \ldots \wedge A_{L} x.
  \end{equation}

  La partie $A_{1}, \ldots, A_{l}$ sera dite partie strictement négative. Et la
  partie $A_{l+1}, \ldots, A_{L}$ sera dite positive.
  \end{proc}

  \section{Application du lemme de saturation ou maximisation à un  MAP non positif}
  Nous allons justifier dans un premier temps le principe de
  saturation ou maximisation sur un exemple
  avec trois variables $x_{1}, x_{2},x_{3}$ et l'atome
  $\{x_{3} \leq a (x_{1} + x_{2})\}$, avec $a < \dun$ (en fait dans cet
  exemple particulier le positionnement de $a$ par rapport à $\dun$
  n'est pas vraiment important). Puis
  nous généraliserons.

  Dans un premier temps nous avons:

  \[
x_{3} \leq a (x_{1} + x_{2}) \Leftrightarrow x \leq A_{1} x
  \]
  avec $x=\left(\bar{c} x_{1} \\ x_{2} \\ x_{3} \ear \right)$ et
  $A_{1}= \left(\bar{ccc} \dun & \dzero & \dzero \\
  \dzero & \dun & \dzero \\
  a & a & \dzero
  \ear \right)$.

  Nous utilisons la décomposition fondamentale suivante:

  $A_{1} x = A x + b$, avec $A=O$ et $b= B u$ où $B=\left(\bar{cc} \dun & \dzero \\
  \dzero & \dun \\
  a & a
  \ear \right)$
  et $u=\left(\bar{c} x_{1} \\ x_{2} \ear \right)$.

  Dans notre cas $A^{*}=I_{3}$. Par application du lemme de
  saturation, il vient: $x = A^{*}B u$. 

  Les variables du vecteur $u$ sont appelées par les automaticiens les
  variables de contrôle. Nous les appelerons ici {\em variables libres} pour
  rappeler qu'elles peuvent prendre n'importe quelle valeur et qu'une fois leur valeur fixée on en déduit la borne supérieure sur $x_{3}$.

  Dans le cas général, nous pouvons définir le concept de variable libre
  par rapport à l'inégalité $x \leq A x$.
  
  \begin{defi}[Variable libre]
    La variable $x_{i}$ est dite libre par rapport au système d'inéquations
    $x \leq A x$ si la ième ligne de $A$ (qui correspond à $x_{i}$) est
    égale à $e_{i}$, le ième vecteur de la base canonique de
    $\real_{\dzero}^{n}$.
  \end{defi}

  Ceci est équivalent à dire que le sommet $i$ (correspondant à la
  variable $x_{i}$) du graphe valué associé à la matrice $A$, $G(A)$,
  n'a pas de prédécesseur à part lui-même.

  Nous généralisons le concept à un ensemble de variables
  libres comme suit.
  Les $k$ premières variables $x_{1}, \ldots, x_{k}$ du vecteur
  $x=(x_{i})_{i=1}^{n}$ sont libres par rapport au système d'inéquations
  $x \leq A x$ si $A$ est de la forme:
  \begin{equation}
    \label{eqmatstructMAPneg}
  A=\left(\bar{cc}
  I_{k} & O_{k,n-k} \\
  B & C
  \ear \right).
  \end{equation}
  Les $k$ premières lignes de $A$ correspondent aux $k$ premiers
  vecteurs de la base canonique de $\real_{\dzero}^{n}$.
  Dans notre cas particulier précédent nous avons: $k=2$,
  $B=\left( \bar{cc} a & a \ear \right)$ et $C=(\dzero)$.

  \subsection{Etape $1$: résolution de $x \leq A_{1} x$}
  \label{subetape1}

  Eliminons déjà deux cas triviaux.
  \bit
\item Cas $1$. Si $A_{1}=I_{n}$ alors par construction toutes
  les lignes des matrices $A_{2}, \ldots, A_{L}$ contiennent au moins un
  terme $\geq \dun$ et le $n$-vecteur colonnes dont les composantes valent
  toutes $\dun$, noté dans ce papier $x^{1}=(\dun)$, est solution évidente.
  Dans l'optique de la recherche
  de toutes les solutions il faut discuter sur les
  matrices $A_{2}, \ldots, A_{L}$. Nous discuterons de ce problème plus en
  détail dans la Section~\ref{secPos}.
  
\item Cas $2$. Il n'y a pas de variables libres dans le système $x \leq A_{1} x$. Alors
  $x=(\dzero)$ est la seule solution. Il s'agit d'une application
  directe du Lemme de saturation avec $b=(\dzero)$ et le cas $4$
  où la matrice $A_{1}$ a toutes ses entrées $< \dun$.
  \eit

  En dehors de ces deux cas on peut supposer qu'il existe
  $1 \leq k < n$ variables libres. Quitte à changer la numérotation des
  varaibles $x_{i}$ on peut supposer que $x_{1}, \ldots , x_{k}$ sont
  libres. Dans ce cas, on a:

  \[
  A_{1} = \left(\bar{cc} I_{k} & O_{k, n-k} \\
  B_{1}     & C_{1}
  \ear \right).
  \]
  Avec $B_{1}$ matrice $(n-k) \times k$ qui a toutes ses
  entrées $ < \dun$. Et, d'après la procédure de
  remplissage~\ref{procRempli}, $C_{1}$ est une matrice $(n-k) \times (n-k)$
  telle que si une ligne $i$ de $B_{1}$ est $\dzero$ alors tous les
  termes de la ligne $i$ dans $C_{1}$ sont $\dzero$ sauf son terme
  diagonal qui vaut $\dun$ (atome $x_{i} \leq x_{i}$). Si la ligne
  $i$ de $B_{1}$ est $\neq \dzero$ alors les termes de la ligne $i$ de
  $C_{1}$ sont $ < \dun$ et le terme diabonal est $=\dzero$. Dans ces
  conditions, il est évident que  $C_{1}$ admet
  une étoile de Kleene.

  Nous avons alors {\bf la décomposition fondamentale}:
\begin{subequations}
  \begin{equation}
x \leq A_{1} x \Leftrightarrow x \leq A x + b,
  \end{equation}
  avec:
  \begin{equation}
  A=\left(\bar{cc} O_{k} & O_{k, n-k} \\
  O_{n-k,k}     & C_{1}
  \ear \right).
  \end{equation}
  et
  \begin{equation}
  b = \left(\bar{c} I_{k} \\ B_{1} \ear\right) u, \; u= \left(\bar{c} x_{1} \\
  \vdots \\ x_{k} \ear \right).
  \end{equation}
  \end{subequations}

  Une application directe du lemme de saturation nous fournit
  l'égalité suivante: $x = A^{*}b$. En remarquant que
  $A^{*} = A=\left(\bar{cc} I_{k} & O_{k, n-k} \\
  O_{n-k,k}     & C_{1}^{*}
  \ear \right)$, il vient:

  \[
x = \left(\bar{c} I_{k} \\ C_{1}^{*} B_{1} \ear\right) u.
  \]
  Notons $J_{1}= C_{1}^{*} B$
  et $T_{1} = \left(\bar{c} I_{k} \\ J_{1} \ear\right)$.

  \begin{rema}
    \label{remImportante}
    Nous avons les équivalences suivantes.
    \bit
  \item[1]. Les variables $u= \left(\bar{c} x_{1} \\ \vdots \\ x_{k} \ear \right)$ sont libres dans le système $x \leq A_{1} x$.
  \item[2]. La matrice $A_{1}$ est de la forme:
    \[
  A_{1} = \left(\bar{cc} I_{k} & O_{k, n-k} \\
  B_{1}     & C_{1}
  \ear \right).
  \]
\item[3]. Les sommets $\{1, \ldots, k\}$ du graphe valué associé à
  la matrice $A_{1}$, $G(A_{1})$, n'ont aucun prédécesseur, à part
  eux-mêmes.
\item[4]. On peut construire une unique matrice $T_{1}$ de la forme
  $\left(\bar{c} I_{k} \\ J_{1} \ear \right)$ associée
  à $A_{1}$ telle que $x= T_{1} u$ est solution de
  $x \leq A_{1} x$, pour tout $u \in \real_{\dzero}^{k}$. Ce qui est équivalent
  à:
  \begin{equation}
    \label{eqvarlibreT}
\forall j=1, \ldots, k: \; t_{1}^{j} \leq A_{1} t_{1}^{j},  
  \end{equation}
  où $t_{1}^{j}$ est le vecteur colonne $j$ de $T_{1}$. 
    \eit
    \end{rema}

  Par construction, nous savons que toutes les matrices de la
  partie négative ont la même forme que la matrice $A_{1}$.

  Pour la partie positive, nous avons deux cas.

  \bit
\item Cas $1$. $A_{l+1}$ a la même forme que $A_{1}$. Et donc, par
  construction toutes les matrices de la partie positive
  ont la même forme que $A_{1}$.
 \item Cas $2$. C'est le cas contraire du cas $1$.
  \eit

  Quitte à permuter certaines matrices dans le système
  $[A_{1}, \ldots, A_{L}]$, on peut supposer qu'il existe
  $l' \geq l$ tel que:
  \[
\forall i=1, \ldots, l':  A_{i} = \left(\bar{cc} I_{k} & O_{k, n-k} \\
  B_{i}     & C_{i}
  \ear \right).
  \]
  avec $C_{i}^{*}$ qui existe (ie pas de circuits de poids $> \dun$).

  En utilisant la même décomposition que pour la matrice
  $A_{1}$, nous construisons une suite de
  matrices $T_{i} = \left(\bar{c} I_{k} \\ J_{i} \ear\right)$
  avec $J_{i}=C_{i}^{*} B_{i}$, $i=1, \ldots , l'$.

  A ce stade, il nous faut construire une matrice $T$ telle que

  \[
\forall i=1, \ldots, l', \; \forall u \in \real_{\dzero}^{k}: T u \leq A_{i} T u.
  \]
  Comme $u$ est un vecteur colonne de variables libres pouvant prendre
  n'importe quelle valeur, il faut construire
  une matrice $T$ telle que:
   \[
\forall i=1, \ldots, l', \;: T \leq A_{i} T. 
  \]
  Comme on doit avoir $x=T u$, la matrice $T$ est de la forme
  $T= \left(\bar{c} I_{k} \\ X \ear \right)$ où $X$ est
  une matrice à calculer.

  Nous avons $\forall i=1, \ldots, l'$:
  \[
  T \leq A_{i} T \Leftrightarrow \left\{ \bar{ccl}
  I_{k} & \leq & I_{k} \\
  X    &  \leq & B_{i} + C_{i} X
  \ear \right.
  \]

  Le même raisonnement que dans la démonstration de lemme de saturation nous
  permet d'obtenir:
  \[
X \leq C_{i}^{*} B_{i} = J_{i}, \; \forall i=1, \ldots, l'.
  \]
  Ce qui veut dire $T \leq T_{i}$, $\forall i = 1, \ldots, l'$.
  Donc, la plus grande matrice $T$ vérifiant ces inégalités est:
  \begin{equation}
T^{\wedge} := T_{1} \wedge \cdots \wedge T_{l'}.
    \end{equation}

  En conclusion, pour tout $u \in \real_{\dzero}^{k}$, $x=T^{\wedge} u$
  est solution de:

  \[
x \leq A_{1} x \wedge \cdots \wedge A_{l'} x.
\]

L'ensemble de toutes les solutions au système
$[A_{1}, \ldots , A_{l'}]$ se déduit de la Proposition suivante.

\begin{propo}
  \label{propvachementImportante}
  L'ensemble $\mathcal{T}_{1,l'}$ de toutes matrices
  $T$ satisfaisant les contraintes $\forall i=1, \ldots, l'$: $T \leq A_{i} T$
  est tel que:
  \begin{equation}
    \label{eqT1lprime}
\mathcal{T}_{1,l'} = \{ T^{\wedge} D, \mbox{ $D$ matrice diagonale $\leq I_{k}$}\}.
    \end{equation}
  
  \end{propo}
\proof
Définissons dès maintenant
\begin{equation}
  \label{eqJwedge}
J^{\wedge}= \bigwedge_{i=1}^{l'}C_{i}^{*} B_{i}.
\end{equation}
Et donc nous avons $T^{\wedge}=\left(\bar{c} I_{k} \\ J^{\wedge} \ear \right)$
qui est la plus grande matrice vérifiant $T^{\wedge} \leq A_{i} T^{\wedge}$,
$\forall i=1, \ldots, l'$. Donc en suivant la même partition que $T^{\wedge}$, toute matrice
$T=\left(\bar{c} D \\ X \ear \right)$ vérifiant les mêmes inégalités
est $\leq T^{\wedge}$. On en déduit déjà que $D \leq I_{k}$. Ceci implique
que $D$ est une matrice diagonale.
 Nous avons $\forall i=1, \ldots, l'$:
  \[
  T \leq A_{i} T \Leftrightarrow \left\{ \bar{ccl}
  D & \leq & D \\
  X    &  \leq & B_{i} D + C_{i} X.
  \ear \right.
  \]
  On en déduit que $X \leq \bigwedge_{i=1}^{l'} (C_{i}^{*} B_{i} D)$.
  Comme $D$ est une matrice diagonale:
  $\bigwedge_{i=1}^{l'} (C_{i}^{*} B_{i} D)= J^{\wedge} D$. Finalement,
  on a
  \[
T \leq T^{\wedge} D,
  \]
  soit encore
  \[
\forall j=1, \ldots, k: t_{j} \leq t_{j}^{\wedge} \; d_{j}.
  \]
  Où $t_{j}$, $t_{j}^{\wedge}$ et $d_{j}$ sont les vecteurs colonnes
  $j$ des matrices $T$, $T^{\wedge}$ et le terme $j$ de la diagonale
  de $D$, respectivement.
  Comme $D$ est une matrice diagonale arbitrire $\leq I_{k}$
  on a: $\{T: T \leq T^{\wedge} D, D \leq I_{k}\} =
  \{T: T = T^{\wedge} D, D \leq I_{k}\}$. Et la Proposition
  est démontrée.
  
\cqfd

Passons à la partie $[A_{l'+1}, \ldots, A_{L}]$, i.e
aux système d'inégalités: $x \leq A_{l'+1} x$ et ...
et $x \leq A_{L} x$. Une
conséquence importante du
résultat de la Proposition~\ref{propvachementImportante} est
que l'on peut raisonner uniquement sur la matrice $T^{\wedge}$.

En écrivant: $x=T^{\wedge} u = t_{1}^{\wedge} x_{1} + \cdots + t_{k}^{\wedge} x_{k}$,
avec $t_{1}^{\wedge}, \ldots , t_{k}^{\wedge}$ les $k$ vecteurs colonne
de $T^{\wedge}$. D'après la Remarque~\ref{remImportante}, en généralisant
le caractère 'libre' d'une variable (cf. (\ref{eqvarlibreT})),
la variable $x_{i}$ est libre (c'est-à-dire qu'elle
peut prendre n'importe quelle valeur)
par rapport au système $[A_{l'+1}, \ldots, A_{L}]$ si et seulement si
\begin{equation}
  \label{xiLibre}
\forall j = l'+1, \ldots, L: t_{i}^{\wedge} \leq A_{j} t_{i}^{\wedge}.
\end{equation}
Remarquons que pour la partie $[A_{1}, \ldots , A_{l'}]$, par construction de
la matrice $T^{\wedge}$, on a déjà:
\[
\forall j = 1, \ldots, l', \forall i=1, \ldots, k: t_{i}^{\wedge} \leq A_{j} t_{i}^{\wedge}.
\]

Nous avons alors deux cas.

\bit
\item Cas $1$. $\forall i=1, \ldots, k$, le système d'inégalités
  (\ref{xiLibre}) n'est pas vérifié. Aucune variable du vecteur 
  $u= \left(\bar{c} x_{1} \\ \vdots \\ x_{k} \ear\right)$ n'est libre dans le
  la partie $[A_{l'+1}, \ldots, A_{L}]$.
  Et donc $u=(\dzero)$ est seule solution possible pour
  la partie $[A_{l'+1}, \ldots, A_{L}]$. Comme $x = T^{\wedge} u$ est
  la plus grande solution (à $u$ fixé) du système $[A_{1}, \ldots , A_{l'}]$
  on a: $x=(\dzero)$, seule solution du système entier:
  $[A_{1}, \ldots , A_{L}]$.

\item Cas $2$. Quitte à renuméroter les variables de $u= \left(\bar{c} x_{1} \\ \vdots \\ x_{k} \ear\right)$, on peut supposer qu'il existe $1 \leq k' \leq k$
  tel que:  $\forall i=1, \ldots, k'$, le système d'inégalités
  (\ref{xiLibre}) est vérifié. 
  
\eit

Nous consacrons la sous-section suivante à l'étude
du deuxième cas.

\subsection{Traitement du Cas $2$: il existe encore $k'$ variables libres}
\label{subCas2}

\noi
a). Si on a $k'=k$, la procédure est finie. En effet, $x = T^{\wedge} u$ est
solution du système entier $[A_{1}, \ldots , A_{L}]$ (équivalent au
MAP initial, $\mathcal{A}$). De la Proposition~\ref{propvachementImportante}
nous déduisons aisément que:

\begin{equation}
\mathcal{A}  =\{ x \in \real_{\dzero}^{n}: x = T^{\wedge} D u, u \in \real_{\dzero}^{k}, D \leq I_{k} \}.
  \end{equation}

\noi
b). Si $1 \leq k' <k$. Définissons alors les trois vecteurs suivants:
$u^{1}= \left(\bar{c} x_{1} \\ \vdots \\ x_{k'} \ear\right)$,
$\overline{u}^{1}= \left(\bar{c} x_{k'+1} \\ \vdots \\ x_{k} \ear\right)$
et $\overline{u}= \left(\bar{c} x_{k+1} \\ \vdots \\ x_{n} \ear\right)$.
Nous pourrions nous contenter de poser $\overline{u}^{1}=(\dzero)$. Donc,
$u= \left(\bar{c} u^{1} \\ \dzero \ear \right)$. Et $x = T^{\wedge} u$.
Mais nous pouvons aller un peu plus loin ici.

En effet, l'étude de la partie $[A_{1}, \ldots , A_{l'}]$ a permis de prouver
que l'ensemble des plus grandes solutions de ce système était $\{T^{\wedge} u, u
\in \real_{\dzero}^{k} \}$, avec $T^{\wedge}= \left(\bar{c} I_{k}
\\ J^{\wedge} \ear \right)$. Où l'on rappelle
que $J^{\wedge}$ est définie par (\ref{eqJwedge}).

Remarquons alors que  cette matrice $J^{\wedge}$ a été
obtenue en posant $\overline{u} = X u$,  $T= \left(\bar{c} I_{k}
\\ X \ear \right)$. Et en résolvant $T \leq A_{i} T$, $i=1, \ldots , l'$. 

Posons:
\[
\overline{u}^{1} = F u^{1},
\]
où $F$ est une $(k-k') \times k'$-matrice à calculer.
Nous avons alors $u = \left(\bar{c} u^{1} \\ \overline{u}^{1} \ear \right)
= T^{1} u^{1}$, avec
\begin{equation}
  \label{eqT1}
  T^{1} = \left(\bar{c} I_{k'} \\ F \ear \right).
  \end{equation}

Nous adoptons le partitionnement de l'espace $\real_{\dzero}^{n}$
suivant les vecteurs $(u^{1}, \overline{u}^{1}, \overline{u})$.

La $n \times k$-matrice $T^{\wedge}$ se réécrit suivant la nouvelle
partition comme suit:

\begin{equation}
  \label{eqpartTwedge}
T^{\wedge}= \left(\bar{cc}
I_{k'} & O_{k', k-k'} \\
O_{k-k',k'} & I_{k-k'} \\
J & K
\ear \right).
\end{equation}

Notons que $\forall F$, la partie $[A_{1}, \ldots , A_{l'}]$ est
toujours vérifiée. 
Soit donc $A$ une $n \times n$-matrice quelconque
de l'ensemble $\{A_{l'+1}, \ldots , A_{L}\}$.
Suivant la partition adoptée nous pouvons écrire:

\[
A = \left(\bar{ccc} A_{1,1} & A_{1,2} & A_{1,3} \\
 A_{2,1} & A_{2,2} & A_{2,3} \\
 A_{3,1} & A_{3,2} & A_{3,3} \\
\ear \right).
\]

Puisque $u^{1}$ est le nouveau vecteur de variables libres
nous exprimons tout en fonction de $u^{1}$.
\[
u = T^{1} u^{1}, \; x = T^{\wedge} u = T^{\wedge} T^{1} u^{1},
\]
avec
\[
T^{\wedge} T^{1} = \left(\bar{c} I_{k'} \\ F \\ J + KF \ear \right),
\]
que l'on réécrit:
\begin{equation}
  \label{eqTwedgeT1}
T^{\wedge} T^{1} = \left(\bar{c} I_{k'} \\ O_{k-k',k'} \\ J\ear \right) +
\left(\bar{c} O_{k'} \\ I_{k-k'} \\ K \ear \right) F.
\end{equation}

A partir de cette expression nous en déduisons l'expression de
$A T^{\wedge} T^{1}$:

\begin{equation}
  \label{eqATwedgeT1}
A T^{\wedge} T^{1} = \left(\bar{c} A_{1,1} + A_{1,3} J \\ A_{2,1} + A_{2,3} J \\
A_{3,1} + A_{3,3} J \ear \right) + \left(\bar{c} A_{1,2} + A_{1,3} K \\
 A_{2,2} + A_{2,3} K \\  A_{3,2} + A_{3,3} K \ear \right) F.
\end{equation}

Comme $x$ est maintenant une fonction du vecteur de variables libres
$u^{1}$ et de la matrice $F$ (à calculer), il nous
faut résoudre:
\[
T^{\wedge} T^{1} \leq A T^{\wedge} T^{1},
\]
ce qui est équivalent à:

\begin{equation}
  \label{eqF1}
 \underbrace{\left(\bar{c} I_{k'} \\ O_{k-k',k'} \\ J\ear \right)}_{u^{1}} +
 \underbrace{\left(\bar{c} O_{k'} \\ I_{k-k'} \\ K \ear \right) F}_{\overline{u}^{1}}
 \leq \underbrace{\left(\bar{c} A_{1,1} + A_{1,3} J \\ A_{2,1} + A_{2,3} J \\
A_{3,1} + A_{3,3} J \ear \right)}_{u^{1}} + \underbrace{\left(\bar{c} A_{1,2} + A_{1,3} K \\
 A_{2,2} + A_{2,3} K \\  A_{3,2} + A_{3,3} K \ear \right) F}_{\overline{u}^{1}}.
\end{equation}

Comme $+$ est l'opérateur $\max$, l'inégalité (\ref{eqF1}) est équivalente
à:

\begin{subequations}
  \begin{equation}
    \label{eqF11}
    \underbrace{\left(\bar{c} I_{k'} \\ O_{k-k',k'} \\ J\ear \right)}_{u^{1}}
    \leq \underbrace{\left(\bar{c} A_{1,1} + A_{1,3} J \\ A_{2,1} + A_{2,3} J \\
A_{3,1} + A_{3,3} J \ear \right)}_{u^{1}} + \underbrace{\left(\bar{c} A_{1,2} + A_{1,3} K \\
 A_{2,2} + A_{2,3} K \\  A_{3,2} + A_{3,3} K \ear \right) F}_{\overline{u}^{1}}.
  \end{equation}
  Et
vv  \begin{equation}
    \label{eqF12}
 \underbrace{\left(\bar{c} O_{k'} \\ I_{k-k'} \\ K \ear \right) F}_{\overline{u}^{1}}
 \leq \underbrace{\left(\bar{c} A_{1,1} + A_{1,3} J \\ A_{2,1} + A_{2,3} J \\
A_{3,1} + A_{3,3} J \ear \right)}_{u^{1}} + \underbrace{\left(\bar{c} A_{1,2} + A_{1,3} K \\
 A_{2,2} + A_{2,3} K \\  A_{3,2} + A_{3,3} K \ear \right) F}_{\overline{u}^{1}}.
    \end{equation}
  \end{subequations}

Puisque $u^{1}$ est le vecteur de variables libres du système $[A_{l'+1}, \ldots,
  A_{L}]$, alors pour toute matrice $A$ de ce système, nous avons:
\[
 \underbrace{\left(\bar{c} I_{k'} \\ O_{k-k',k'} \\ J\ear \right)}_{u^{1}}
    \leq \underbrace{\left(\bar{c} A_{1,1} + A_{1,3} J \\ A_{2,1} + A_{2,3} J \\
A_{3,1} + A_{3,3} J \ear \right)}_{u^{1}},
    \]
    et l'inégalité (\ref{eqF11}) est donc vérifiée $\forall F$.

    En notant
    \begin{equation}
      \label{defZ}
      Z= \left(\bar{c} I_{k-k'} \\ K \ear \right) F,
    \end{equation}
    
    \begin{equation}
      \label{defB}
      B= \left(\bar{c} A_{2,1} + A_{2,3} J \\
      A_{3,1} + A_{3,3} J \ear \right)
    \end{equation}
    et en remarquant que:
    \[
    \left(\bar{c} A_{2,2} + A_{2,3} K \\  A_{3,2} + A_{3,3} K \ear \right) F =
    C \left(\bar{c} I_{k-k'} \\ K \ear \right) F = C Z, 
    \]
    où $C$ est la matrice carrée d'ordre $n-k'$ suivante:
    \begin{equation}
      \label{defC}
C= \left(\bar{cc} A_{2,2} & A_{2,3} \\ A_{3,2} & A_{3,3} \ear \right). 
\end{equation}

Il vient que l'inégalité (\ref{eqF12}) est équivalente à:

\begin{subequations}
  \begin{equation}
    \label{eqF121}
Z \leq C Z + B
\end{equation}
et
\begin{equation}
  \label{eqF122}
\left(\bar{c} I_{k-k'} \\ K \ear \right) F = Z.
\end{equation}
\end{subequations}

Nous résolvons ce problème (\ref{eqF121})-(\ref{eqF122}) comme suit.

\noi
a/ Résolution de (\ref{eqF121}). Nous avons deux cas.

\bit
\item Cas $1$. $C^{*}$ existe (ie dans le graphe valué associé à $C$, $G(C)$
  il n'y a pas de circuit de poirds $> \dun$). Dans ce cas
  la plus grande solution est: $Z = C^{*} B$
\item Cas $2$. $C^{*}$ n'existe pas. Prendre $Z=B$.   

\eit

\noi
b/ Résolution de (\ref{eqF122}) connaissant $Z$.
La théorie de la résiduation nous fournit au moins la plus grande matrice
$F$ (au sens de l'ordre partiel $\leq$ sur
les matrices) sous-solution, ie la matrice $F$ qui vérifie
\begin{equation}
\label{eqResidF}
\left(\bar{c} I_{k-k'} \\ K \ear \right) F \leq Z
\end{equation}
que l'on notera:
\begin{equation}
  \label{eqFA}
F(A) = \left(\bar{c} I_{k-k'} \\ K \ear \right) \backslash Z.
\end{equation}

Par le même procédé, nous pouvons construire les matrices
$F(A_{i})$, $i=l'+1, \ldots , L$. En remarquant que toute matrice $F \leq F(A)$
par définition de l'opérateur résiduation vérifie encore l'inéquation
(\ref{eqResidF}), la plus grande matrice, notée $F^{\wedge}$, qui vérifie
toutes les contraintes est définie par:

\begin{equation}
  \label{eqFwedge}
F^{\wedge}= F(A_{l'+1}) \wedge \ldots \wedge F(A_{L}).
\end{equation}

Finalement, nous avons les plus grandes solutions du MAP
$\mathcal{A}$:

\begin{equation}
  \label{eqsolgen}
\sup_{\leq}\{x \in \real_{\dzero}^{n}: \mbox{$\mathcal{A}$ est vérifié} \}= \{
x= \left(\bar{c} I_{k'} \\ F^{\wedge} \\ J + KF^{\wedge} \ear \right) u^{1}, \;
u^{1} \in \real_{\dzero}^{k'} \}.
\end{equation}

En utilisant la relation (\ref{eqTwedgeT1}, avec $F=F^{\wedge}$):
\[
\left(\bar{c} I_{k'} \\ F^{\wedge} \\ J + KF^{\wedge} \ear \right) =
\left(\bar{c} I_{k'} \\ O_{k-k',k'} \\ J \ear \right) + \left(\bar{c} O_{k'} \\
I_{k-k'} \\ K \ear \right) F^{\wedge}
\]

Nous remarquons que la matrice $\left(\bar{c} I_{k'} \\ O_{k-k',k'}
\\ J \ear \right)$ est une sous-matrice correspondant aux $k'$
premières colonnes de $T^{\wedge}$ (cf. (\ref{eqpartTwedge})).
Le résultat de la Proposition~\ref{propvachementImportante}
s'applique à cette sous-matrice. Et d'autre part, par définition
de la matrice $F^{\wedge}$ nous déduisons toutes les solutions du
MAP $\mathcal{A}$ et nous avons le résultat final suivant.

\begin{equation}
  \label{eqFIN}
\mathcal{A} = \{ x \in \real_{\dzero}^{n}: x = T^{\wedge} \left(\bar{c} D \\ F \ear \right) u^{1}, u^{1} \in \real_{\dzero}^{k'}, D \leq I_{k'}, F \leq F^{\wedge} \}
\end{equation}
en rappelant ici le partitionnement adopté pour $T^{\wedge}$
(cf. (\ref{eqpartTwedge})):

\[
T^{\wedge}= \left(\bar{cc}
I_{k'} & O_{k', k-k'} \\
O_{k-k',k'} & I_{k-k'} \\
J & K
\ear \right).
\]

\section{Cas d'un MAP positif}
\label{secPos}
Dans cette section nous considérons un MAP $\mathcal{A}$ sur $\real_{\dzero}^{n}$
dit positif. C'est-à-dire que par rapport à la procédure de
remplissage des matrices associées à $\mathcal{A}$ (cf. sous-section~\ref{subCasGenAlgebric})
nous avons $A_{1}=I_{n}$ et toutes les matrices
$A_{2}, \ldots , A_{L}$ du sytème équivalent à $\mathcal{A}$ ont tous
les vecteurs-lignes qui les composent avec au moins une composante $\geq \dun$.
En conséquence, en prenant $x^{1}=(\dun)$ nous avons:
\[
\forall i=1, \ldots, L: x^{1} \leq A_{i} x^{1}.
\]

La première remarque qui vient est de se rendre compte que la structure matricielle 
imposée (\ref{eqmatstructMAPneg}) pour avoir des solutions non
triviales $x$, $ x \neq (\dzero)$ n'a plus lieu d'être ici. De toute
manière même si le MAP positif a des variables libres et que la structure
matricielle (\ref{eqmatstructMAPneg}) existe il faudrait que la matrice $C^{*}$ existe, i.e.
son graphe associé $G(C)$ devra donc être un DAG. Ce qui n'est pas toujours possible.


\subsection{Analyse de la difficulté: les circuits élémentaires}
\label{ideediff}
La recherche de toutes les solutions d'un MAP positif est liée aux
{\bf circuits élémentaires} des graphes valués associés aux matrices du MAP.

En effet, considérons sur $\real_{\dzero}^{5}$ le système $x \leq A x$ avec:
\begin{equation}
  \label{exmatpos}
  A=\left(\bar{ccccc}
  \dzero & \dun & \dzero & \dzero & \dzero \\
  \dzero & \dzero & \dun & \dzero & \dzero \\
  \dun & \dun & \dzero & \dzero & \dzero \\
  \dzero & \dzero & \dun & \dzero & \dun \\
  \dzero & \dzero & \dzero & \dun & \dzero
  \ear \right)
\end{equation}

Le graphe $G(A)$ associé à la matrice $A$
possède deux composantes fortement connexes
$CFC(1) = \{1,2,3\}$ et $CFC(2)= \{4,5\}$. Le graphe orienté sans cycles
(ou DAG) des CFC est: $CFC(1) \longrightarrow CFC(2)$.
Nous remarquons que la $CFC(1)$ possède le {\bf circuit élémentaire} $\{2, 3\}$.
Or, si l'on prend le vecteur $x = \left(\bar{c} \dzero \\ \dun \\ \dun \\ \dzero \\ \dzero
\ear \right)$ dont le support est le circuit $\{2, 3\}$. On a
$A x = \left(\bar{c} \dun \\ \dun \\ \dun \\ \dun \\ \dzero
\ear \right)$. Et donc: $x \leq A x$. 

\subsection{Une méthode pour trouver des solutions non triviales}
\label{subSolnontriviales}
Dans un premier temps nous traitons un seul système d'inégalités $x \leq Ax$, où $A =A_{2}$.

La recherche de toutes les solutions d'un MAP
positif repose sur l'énumération de
des graphes associés aux matrices de ce MAP.
Et ce nombre n'est donc pas forcément polynomial en le nombre de variables
et le nombre d'arcs dans les graphes associés aux matrices du MAP.

Nous proposons alors une approche pour tenter de fournir des solutions
$x$ non triviales, i.e. $\neq x^{1}$. Elle repose sur la notion de matrice
monomiale associée naturellement à la notion de circuit
élémentaire. Nous renvoyons e.g. à \cite[pp. 45-46]{kn:Gaubertthese} pour la
définition suivante.

\begin{defi}[Matrice monomiale]
  $A$ est une matrice carrée monomiale si elle s'écrit sous la forme
  produit suivante:
  \begin{equation}
A = D P
  \end{equation}
  avec
  $D$ une matrice diagonale inversible (i.e. tous ses termes
  diagonaux sont $\neq \dzero$), et $P$ une matrice de permutation. 
\end{defi}

Il est important de noter que toute matrice $A$ monomiale possède une unique inverse, notée
$A^{-}$, qui vérifie l'égalité bien connue en algèbre linéaire dont la
version maxplusienne est $A A^{-} = A^{-} A = I_{n}$. Le matrice $A^{-}=(a_{i,j}^{-}$ étant
définie par:

\begin{equation}
  \label{defA-}
  \forall i,j \in \{1, \ldots , n\}, \; a_{i,j}^{-} = \left\{ \bar{cc} a_{j,i}^{-1} & \mbox{si $a_{j,i} \neq \dzero$} \\
                                                                     \dzero & \mbox{si $a_{j,i} = \dzero$}
  \ear \right.
  \end{equation}

Remarquons également que nous avons clairement l'équivalence logique suivante.

\begin{propo}
  La matrice $A$ est monomiale si et seulement si
  son graphe $G(A)$ ne possède qu'un seul circuit
  élémentaire avec éventuellement des variables sans prédécesseurs (i.e. des variables
  libres).
\end{propo}

De ces deux précédents résultats nous avons alors le résultat remarquable suivant.

\begin{theo}
  L'ensemble de toutes les solutions du cône $\{x \in \real_{\dzero}^{*}: x \leq Ax\}$ est l'ensemble
  des combinaisons $(max,+)$-linéaires des colonnes
  de $(A^{-})^{*}$, i.e. $(A^{-})^{*} \real_{\dzero}^{n}$.
  \end{theo}
\proof Comme les fonctions $x \mapsto Ax$ et $x \mapsto A^{-} x$ sont
croissantes nous avons l'équivalence
logique suivante:

\begin{equation}
x \leq A x \Leftrightarrow A^{-} x \leq x.
  \end{equation}

Or, la matrice $A^{-}$ n'a que des termes $\leq \dun$. Elle ne possède donc pas de circuits de
poids $ > \dun$ et donc $(A^{-})^{*}$ existe. Un résultat classique en algèbre $(\max,+)$ nous permet alors
de conclure que toutes les solutions du sytème
$A^{-} x \leq x$ s'écrivent: $x = (A^{-})^{*} y$, $y \in \real_{\dzero}^{n}$. Et la preuve est
achevée. \cqfd

Dans le cas où la matrice $A$ n'est pas monomiale nous proposons l'approche suivante.

\bit
\item[1]. Construire la matrice $A^{-}$ définie par (\ref{defA-}). Par construction, 
  dans le cas général nous avons:
  \begin{equation}
I \leq A^{-} A.
    \end{equation}
\item Au lieu de traiter le système $x \leq A x$ nous traitons le système
  $A^{-} x \leq x$. Et donc, il nous faut les étapes suivantes.
  \bit
\item[(2.a)] Calculer la matrice $(A^{-})^{*}$ que l'on notera $A^{-*}$

\item[(2.b)]. Notons $a_{i}$ le ième vecteur ligne de $A$. Et notons
$a_{j}^{-*}$ le jème vecteur colonne de $A^{-*}$. Puis définissons
 le $n$-vecteur colonne $a_{j}^{\#}$ en notant 

  \begin{equation}
    \forall i \in \{1, \ldots , n\}, \; a_{i,j}^{\#} = \left\{ \bar{cc} a_{i,j}^{-*} & \mbox{si $a_{i} a_{j}^{-*} \leq \dun$} \\
    \dzero    & \mbox{sinon.}
    \ear \right.
    \end{equation}

  Ce vecteur vérifie alors
  \begin{equation}
    \label{yo!}
    a_{j}^{\#} \leq A a_{j}^{\#}.
    \end{equation}

  Au total, à partir des $n$ vecteurs colonnes de
  $A^{-*}$ nous générons $n$ vecteurs $a_{j}^{\#}$ que nous stockons
  dans la matrice notée $A^{\#}$.
  \eit
\eit

Cette méthode présente l'avantage d'être facilement généralisabe au cas
$L \geq 3$ comme suit. A chaque système $x \leq A_{k}x$ on associe
le système $ A_{k}^{-} x \leq x$ où $A_{k}^{-}$ est définie par (\ref{defA-},
$A \leftrightarrow A_{k}$), $k \in \{1, \ldots , L\}$. Il suffit ensuite de remarquer
que nous avons l'équivalence logique suivante:

\begin{subequations}
\begin{equation}
\forall k \in \{1, \ldots ,L \}, \; A_{k}^{-} x \leq x \Leftrightarrow A^{-}x \leq x,
  \end{equation}
où
\begin{equation}
  \label{autreA-}
A^{-}:= A_{1}^{-} + \cdots + A_{L}^{-}.
  \end{equation}
\end{subequations}

Il suffit alors d'appliquer la même méthode mais cette fois avec
la matrice $A^{-}$ définie par (\ref{autreA-}) avec l'étape (2.b)
remplacée par l'étape suivante:

\bit
\item[(2.b')]. Pour tout $k \in \{1, \ldots, L\}$ notons
  $a_{k:i}$ le vecteur ligne $i$ de la matrice $A_{k}$. Et notons
  comme précédemment $a_{j}^{-*}$ le jème vecteur colonne de $A^{-*}$.
   Puis définissons
 le $n$-vecteur colonne $a_{j}^{\#}$ en notant 

  \begin{equation}
    \forall i \in \{1, \ldots , n\}, \; a_{i,j}^{\#} = \left\{ \bar{cl} a_{i,j}^{-*} & \mbox{si $\forall k \in
      \{1, \ldots , L\}$, $a_{k:i} a_{j}^{-*} \leq \dun$} \\
    \dzero    & \mbox{sinon.}
    \ear \right.
    \end{equation}

  Ce vecteur vérifie alors
  \begin{equation}
    \label{yoda!}
\forall k \in \{1, \ldots, L\}: \;    a_{j}^{\#} \leq A_{k} a_{j}^{\#}.
    \end{equation}

\eit

\section{Analyse de la complexité}
\label{secComplexity}
Rappelons que le MAP est constitué de $m$ atomes et $n$ varibles
$x_{1}, \ldots, x_{n}$. Pour chaque variable $x_{i}$ nous avons
$m_{i}^{-}$ atomes strictement négatifs et $m_{i}^{+}$ atomes
positifs. Nous avons l'égalité suivante:
$\sum_{i=1}^{n} (m_{i}^{-} + m_{i}^{+}) = m$. Notons
$\overline{m}= \max_{i=1}^{n} m_{i}^{-}$. Nous avons alors deux cas.

\bit
\item Cas $1$: $\overline{m} > 0$. La plus grosse complexité
  provient alors de la décomposition fondamentale et du calcul
  de l'étoile de Kleene d'une matrice carrée de taille $n$ au maximum. En
  effet le reste de la méthode consiste à calculer des $\min$
  de matrices ou vecteurs et de faire des comparaisons termes à termes.
  Dans la sous-section~\ref{subCas2}, il nous faut détecter si la matrice
  $C$ défine par (\ref{defC}) n'a pas de circuit de poids $ > \dun$.
  L'algorithme eg de Bellman-Ford permet de résoudre ce problème
  en $\mathcal{O}(n^{3})$. 

  Donc, la recherche de toutes les solutions se fait en
  $\mathcal{O}( \overline{m} n^{3})$.

\item Cas $2$: $\overline{m}=0$. Dans ce cas tous les atomes sont positifs. Si l'on se
  contente d'une solution non nulle alors $x^{1}=(\dun)$ est toujours
  solution. La recherche de toutes les solutions est liée à l'énumération
  des circuits élémentaires dans le graphe associé aux matrices du MAP
  (cf. sous-section~\ref{ideediff}). Cela étant nous avons proposé
  une méthode (cf. sous-section~\ref{subSolnontriviales}) qui évite
  l'énumération des circuits élémentaires. Elle est fondée sur l'association naturelle
  $\mbox{circuit élémentaire} \leftrightarrow \mbox{matrice monomiale}$.
  Cette méthode ne peut fournir qu'au plus $n$ solutions
  non triviales, i.e. $\neq x^{1}$, avec une complexité polynomiale en
  $O(n^{3})$ qui correspond au calcul d'une étoile de Kleene d'une
  $(n \times n)$-matrice.

\eit

\section{Retour sur l'exemple de base}
\label{secEDB}

Cette section est organisée comme suit. 
Dans les sous-sections~\ref{subttSprime} et ~\ref{systemeS''} Nous appliquons
la méthode proposée sur l'exemple de
la sous-section~\ref{subEDB}. Et dans la sous-section~\ref{subMappositif}
nous montrons comment calculer des solutions non triviales d'un MAP positif
en utilisant la méthode développée en sous-section~\ref{subSolnontriviales}.

\subsection{Traitement du système $S'$, sous-section~\ref{subEDB}}
\label{subttSprime}
Le MAP $S'$ est $S'= \{x_{3} \leq (-10) + x_{1}, x_{4} \leq (-1) + \max(x_{2},x_{3}), x_{2} \leq x_{4}, 
x_{4} \leq 25 + \max(x_{2},x_{3}), x_{1} \leq 9+ \max(x_{2}, x_{3}) \}$.

En utilisant la procédure de remplissage des matrices de
la sous-section~\ref{subCasGenAlgebric} liée à la constitution des
listes des atomes strictement négatifs et positifs, on obtient:

\[
x \leq A_{1} x \wedge A_{2} x,
\]
avec $x= \left(\bar{c} x_{1} \\ \vdots \\ x_{4} \ear \right)$. Pour la
partie strictement négative, nous avons:
\[
A_{1}= \left(\bar{cc} I_{2} & O_{2,2} \\
B_{1} & C_{1} \ear \right), \; B_{1}= \left(\bar{cc} 10^{-1} & \dzero \\
\dzero & 1^{-1} \ear\right), \; C_{1}= \left(\bar{cc} \dzero & \dzero \\
1^{-1} & \dzero \ear \right).
\]
Pour la partie positive, nous avons:

\[
A_{2} = \left(\bar{cccc} \dzero & 9 & 9 & \dzero \\
\dzero & \dzero & \dzero & \dun \\
\dzero & \dzero & \dun & \dzero \\
\dzero & 25 & 25 & \dzero
\ear \right).
\]

D'après la forme de la matrice $A_{1}$ le vecteur $u=\left(\bar{c} x_{1} \\ x_{2}
\ear \right)$ est un vecteur de variables libres. En développant
$C_{1}^{*} B_{1}= \left(\bar{cc} 10^{-1} & \dzero \\ 11^{-1} & 1^{-1} \ear \right)$
nous avons
\[
\sup_{\leq}\{x: x \leq A_{1} x \} = \{x= T^{\wedge} u, u \in \real_{\dzero}^{2} \},
\]
avec
\begin{equation}
  \label{defTwedgeEx}
T^{\wedge}= \left(\bar{cc} \dun & \dzero \\
\dzero & \dun \\
 10^{-1} & \dzero \\ 11^{-1} & 1^{-1} \ear \right).
\end{equation}

\noi
$x_{1}$ est-elle libre pour le système $[A_{2}]$ ? \\
A t-on: $t_{1}^{\wedge} = \left(\bar{c} \dun \\ \dzero \\ 10^{-1} \\ 11^{-1}
\ear \right) \leq
A_{2} t_{1}^{\wedge} = \left(\bar{c} 1^{-1} \\ 11^{-1} \\ 10^{-1} \\ 15 \ear \right)$ ? réponse non. \\

\noi
$x_{2}$ est-elle libre pour le système $[A_{2}]$ ? \\
A t-on: $t_{2}^{\wedge} = \left(\bar{c} \dzero\\ \dun \\ \dzero  \\ 1^{-1} \ear
\right) \leq
A_{2} t_{2}^{\wedge} = \left(\bar{c} 9 \\ 1^{-1} \\ \dzero \\ 25 \ear \right)$ ?
réponse non. \\

En conclusion, $u=(\dzero)$ et donc $x=(\dzero)$ seule solution de
$[A_{1}, A_{2}]$, ie de $S'$.

\subsection{Traitement du système $S''$, sous-section~\ref{subEDB}}
\label{systemeS''}
Dans ce système seule la partie positive change. Et nous avons la
matrice $A_{2}$ remplacée par $A_{2}^{'}$, avec
\[
A_{2}^{'}= \left(\bar{cccc} \dzero & 12 & 12 & \dzero \\
\dzero & \dzero & \dzero & \dun \\
\dzero & \dzero & \dun & \dzero \\
\dzero & 25 & 25 & \dzero
\ear \right).
\]

Nous nous posons à nouveau les questions suivantes.

\noi
$x_{1}$ est-elle libre pour le système $[A_{2}^{'}]$ ? \\
A t-on: $t_{1}^{\wedge} = \left(\bar{c} \dun \\ \dzero \\ 10^{-1} \\ 11^{-1}
\ear \right) \leq
A_{2}^{'} t_{1}^{\wedge} = \left(\bar{c} 2 \\ 11^{-1} \\ 10^{-1} \\ 15 \ear \right)$ ? réponse oui. \\

\noi
$x_{2}$ est-elle libre pour le système $[A_{2}^{'}]$ ? \\
A t-on: $t_{2}^{\wedge} = \left(\bar{c} \dzero\\ \dun \\ \dzero  \\ 1^{-1} \ear
\right) \leq
A_{2}^{'} t_{2}^{\wedge} = \left(\bar{c} 12 \\ 1^{-1} \\ \dzero \\ 25 \ear \right)$ ?
réponse non. \\

Nous pouvons nous contenter de poser
$u= \left(\bar{c} x_{1} \\ \dzero \ear \right)$. Dans ce cas, nous
avons: $x= t_{1}^{\wedge} x_{1} =
\left(\bar{c} \dun \\ \dzero \\ 10^{-1} \\ 11^{-1} \ear \right) x_{1}$
qui est bien solution du système $[A_{1}, A_{2}^{'}]$ pour tout $x_{1}
\in \real_{\dzero}$.

Mais nous pouvons aussi appliquer les
résultats et la plupart des notations de la sous-section~\ref{subCas2}.
La nouvelle partition de l'espace des coordonnées $\real_{\dzero}^{4}$ est:
$u^{1}=(x_{1})$, $\overline{u}^{1}=(x_{2})$ et $\overline{u}=\left(\bar{c} x_{3} \\ x_{4} \ear \right)$.
Et nous posons $\overline{u}^{1} = F u^{1}$ où $F =(f)$ est à calculer.

Pour cette nouvelle partition, $T^{\wedge}$ se réécrit:

\[
T^{\wedge}=\left(\bar{cc} \dun & \dzero \\ \dzero & \dun \\ J & K \ear \right),
\; J=\left(\bar{c} 10^{-1} \\ 11^{-1} \ear \right), \; K= \left(\bar{c} \dzero \\ 1^{-1} \ear \right).
\]

De même la matrice $A_{2}^{'}$ se réécrit:
\[
A_{2}^{'} = \left(\bar{ccc} A_{1,1} & A_{1,2} & A_{1,3} \\
 A_{2,1} & A_{2,2} & A_{2,3} \\
 A_{3,1} & A_{3,2} & A_{3,3} \\
\ear \right),
\]
avec: \\
$A_{1,1}=\left(\bar{c} \dzero \ear \right)$, $A_{1,2}=\left(\bar{c} 12 \ear \right)$, $A_{1,3}=\left(\bar{cc} 12 & \dzero \ear \right)$, 
$A_{2,1}=\left(\bar{c} \dzero \ear \right)$, $A_{2,2}=\left(\bar{c} \dzero \ear \right)$, $A_{2,3}=\left(\bar{cc} \dzero & \dun \ear \right)$, 
$A_{3,1}=\left(\bar{c} \dzero \\ \dzero \ear \right)$, $A_{3,1}=\left(\bar{c} \dzero \\ 25 \ear \right)$, $A_{3,3}=\left(\bar{cc} \dun & \dzero \\ 25 & \dzero \ear \right)$.

Nous devons alors résoudre le système (\ref{eqF121})-(\ref{eqF122})
avec ici:
\bit
\item $F=(f)$, $f$ scalaire inconnu.
\item $Z = \left(\bar{c} \dun \\ \dzero \\ 1^{-1} \ear \right) F$ (cf. (\ref{defZ}))

\item $B = \left(\bar{c} 11^{-1} \\ 10^{-1} \\ 15 \ear \right)$ (cf. (\ref{defB}))
\item $C= \left(\bar{ccc}  \dzero & \dzero & \dun \\
  \dzero & \dun & \dzero \\
  25 & 25 & \dzero
  \ear \right)$ (cf. (\ref{defC})).
  
 \eit 

 Nous voyons immédiatement que $C^{*}$ n'existe pas. Nous sommes
 dans le Cas $2$ de la résolution de (\ref{eqF121}). Et nous posons
 \[
Z = B. 
\]

\subsubsection{Application de la théorie de la résiduation}
\label{subsubtheoResiduation}

Le calcul de $F(A_{2}^{'}) = (f(A_{2}^{'}))$ définie par
(\ref{eqFA}) s'effectue en appliquant
la formule fournie par la théorie de la Résiduation que nous détaillons
ci-après. Ici, nous avons juste à calculer le scalaire
$f(A_{2}^{'})$ défini par:

\[
f(A_{2}^{'}) = \left(\bar{c} \dun \\ \dzero \\ 1^{-1} \ear \right) \backslash
\underbrace{ \left(\bar{c} 11^{-1} \\ 10^{-1} \\ 15 \ear \right)}_{B}.
\]

Le vecteur colonne $\theta=\left(\bar{c} \dun \\ \dzero \\ 1^{-1} \ear \right)$
est transformé en le vecteur ligne suivant:
$\theta^{-\top}=\left(\bar{ccc} \dun & \infty & 1 \ear \right)$.
Puis nous effectuons le produit $\theta^{-\top}$ par le vecteur $B$
noté $\theta^{-\top} \otimes' B$ dans l'algèbre $(\min, +)$ mais en utilisant
les notations de l'algèbre usuelle:
\[
\theta^{-\top} \otimes' B =\min( 0 + (-11), \infty + (-10), 1 + 15)= -11. 
\]

Avec les notations maxplusiennes de cet article nous avons:
$f^{\wedge}=f(A_{2}^{'})= 11^{-1}$ (voir définition de $f^{\wedge}$
(\ref{eqFwedge})).

\subsubsection{Synthèse pour le système $S''$}

Le système matricel $[A_{1},A_{2}^{'}]$ est équivalent au MAP $S''$.
Donc, en application du résultat (\ref{eqsolgen}), schant que:
\[
T^{\wedge} \left(\bar{c} \dun \\ f^{\wedge}\ear \right)= \left(\bar{c}
\dun \\ 11^{-1} \\ 10^{-1} \\ 11^{-1} \ear \right),
\]
nous avons
\[
\sup_{\leq}\{x \in \real_{\dzero}^{4}: \mbox{$S''$ est vérifié} \}= \{\left(\bar{c}
\dun \\ 11^{-1} \\ 10^{-1} \\ 11^{-1} \ear \right) \; x_{1}, \; x_{1} \in \real_{\dzero} \}.
\]

Toutes les solutions de $S''$ sont alors obtenues en appliquant le
résultat (\ref{eqFIN}), ainsi il vient:
\begin{equation}
S'' = \{x \in \real_{\dzero}^{4}; x= T^{\wedge} \left(\bar{c} d \\ f \ear \right) x_{1}, x_{1} \in \real_{\dzero}, d \leq \dun , f \leq f^{\wedge} \}, 
\end{equation}
en rappelant que la matrice $T^{\wedge}$ est définie
par (cf. (\ref{defTwedgeEx})):
\[
T^{\wedge}= \left(\bar{cc} \dun & \dzero \\
\dzero & \dun \\
 10^{-1} & \dzero \\ 11^{-1} & 1^{-1} \ear \right).
\]

\subsection{Exemple d'un MAP positif}
\label{subMappositif}

Dans cette sous-section, nous déroulons notre méthode présentée en sous-section~\ref{subSolnontriviales}
pour trouver des solutions
non triviales au MAP $S'$, sous-section~\ref{subEDB}, oté de sa partie négative.
Autrement dit, nous considérons le MAP positif $[A_{2}]$ en rappelant que
la matrice $A_{2}$ est la matrice suivante:

\[
A_{2} = \left(\bar{cccc} \dzero & 9 & 9 & \dzero \\
\dzero & \dzero & \dzero & \dun \\
\dzero & \dzero & \dun & \dzero \\
\dzero & 25 & 25 & \dzero
\ear \right).
\]
La matrice $A_{2}^{-}$ est définie par (\ref{defA-}, avec $n=4$, $A \leftrightarrow A_{2}$), i.e. :

\[
A_{2}^{-} = \left(\bar{cccc} \dzero & \dzero & \dzero & \dzero \\
9^{-1} & \dzero & \dzero & 25^{-1} \\
9^{-1} & \dzero & \dun & 25^{-1} \\
\dzero & \dun & \dzero & \dzero 
\ear \right).
\]

Le calcul de l'étoile de Kleene de $A_{2}^{-}$, notée $A_{2}^{-*}$, nous fournit
alors:

\[
A_{2}^{-*} = \left(\bar{cccc} \dun & \dzero & \dzero & \dzero \\
9^{-1} & \dun & \dzero & 25^{-1} \\
9^{-1} & 25^{-1} & \dun & 25^{-1} \\
9^{-1} & \dun & \dzero & \dun 
\ear \right).
\]

Le lecteur vérifiera aisément que $A_{2}^{-*}  \leq A_{2} A_{2}^{-*}$. Et donc la matrice
$A_{2}^{\#}$ est égale à la matrice $A_{2}^{-*}$. Ainsi toutes les
combinaisons $(\max,+)$-linéaires $x$ des colonnes
de $A_{2}^{\#}$ seront solution de $x \leq A_{2} x$. Enfin, remarquons que la
somme des colonnes de $A_{2}^{\#}$ fournit la solution triviale $x^{1}=(\dun)$. 

\bibliographystyle{plain}
\bibliography{ref_lt}

\begin{thebibliography}{1}

\bibitem{kn:Bac-cooq}
F.~Baccelli, G.~Cohen, G.J. Olsder, and J-P. Quadrat.
\newblock {\em {Synchronization and Linearity}}.
\newblock John Wiley and Sons, 1992.

\bibitem{kn:Bezemetal2008}
M.~Bezem, R.~Nieuwenhuis, and E.~Rodriguez-Carbonell.
\newblock {The Max-Atom Problem and Its Relevance}.
\newblock {\em LPAR08}, 2008.
\newblock (47-61).

\bibitem{kn:Gaubertthese}
S.~Gaubert.
\newblock {\em {Th\'eorie des Syst\`emes Lin\'eaires dans les dio\"ides.}}
\newblock PhD thesis, Ecole des Mines de Paris, 1992.

\end{thebibliography}

\end{document}